\documentclass[12pt]{amsart}

\newtheorem{theorem}{Theorem}[section]
\newtheorem{lemma}[theorem]{Lemma}
\newtheorem{prop}[theorem]{Proposition}
\newtheorem{corollary}[theorem]{Corollary}
\theoremstyle{definition}   
\newtheorem{definition}{Definition}
\newtheorem{example}[theorem]{Example}

\theoremstyle{remark}

\numberwithin{equation}{section}

\usepackage{times}
\usepackage{enumerate}
\usepackage{graphicx,adjustbox}
\usepackage{hyperref}
\usepackage{amsmath,amssymb}
\usepackage{amscd}
\usepackage{graphicx}
\usepackage[all]{xy}

\title[Derivation modules for Sums and Gluings]
{Derivation modules for Sum and Gluing}
\author{
Joydip Saha
\and
Indranath Sengupta
}
\date{}

\address{\small \rm  Discipline of Mathematics, IIT Gandhinagar, Palaj, Gandhinagar, 
Gujarat 382355, INDIA.}
\email{saha.joydip56@gmail.com}
\thanks{The first author thanks SERB, Government of India for the Research Associate 
position at IIT Gandhinagar, through the research project EMR/2015/000776.}

\address{\small \rm  Discipline of Mathematics, IIT Gandhinagar, Palaj, Gandhinagar, 
Gujarat 382355, INDIA.}
\email{indranathsg@iitgn.ac.in}
\thanks{The second author is the corresponding author. This research is supported 
by the research project EMR/2015/000776, sponsored by SERB.}

\date{}



\begin{document}

\begin{abstract}
In this paper we explicitly compute the derivation module of quotients of polynomial rings by 
ideals formed by the sum or by some other gluing technique. We discuss cases of monomial 
ideals and binomial ideals separately. 
\end{abstract}

\maketitle
Let $k$ denote a field and $R=k[x_{1},x_{2},\ldots,x_{n}]$ the polynomial ring over $k$. 
Our aim in this paper is to study the structure of the Derivation module for three 
interesting situations: (i) sum of ideals which are intersecting transversally, 
(ii) the join of simplicial complexes, (iii) the gluing of numerical semigroups. 
This study is an attempt to understand and answer the following question: Given ideals 
$I$ and $J$ in $R$, what can one say about the structure 
of the derivation module $\mbox{Der}_{k}(R/I+J, R/I+J)$? This is difficult to answer 
in general and in this paper we aim to answer this question partially for some good 
cases, making use of existing results by Brumatti and Simis in \cite{bs} and by Kraft 
in \cite{kraft}. In the first section we use the result proved by Brumatti and Simis 
\ref{brusim} (from  \cite{bs}) to study $\mbox{Der}_{k}(R/I+J, R/I+J)$, where $I$ 
and $J$ are monomial ideals with the property that $I\cap J = IJ$. As a consequence, 
we derive the structure of the derivation module of the Stanley-Reisner ring of the 
simplicial complex obtained by the join of two simplicial complexes. The second section 
is devoted to the study of the structure of the derivation module of the glued numerical 
semigroup ring formed by the gluing of numerical semigroups. We use a theorem of Kraft \cite{kraft} to prove our result.
\medskip

\section{Sum of monomial ideals}
Ideals $I$ and $J$ 
in $R$ are said to intersect transversally if $I\cap J = IJ$. Transversal intersection of monomial ideals 
and ideals of the form $I+J$, where $I$ and $J$ intersect transversally have been studied in detail in 
\cite{ssttrans} and \cite{sstsum} respectively. Syzygies of ideals of the form $I+J$ can be read easily 
when $I$ and $J$ intersect transversally. This section is devoted to describe the structure of the derivation 
module $\mbox{Der}_{k}(R/I+J, R/I+J)$, where $I$ and $J$ are monomial ideals with the property that 
$I\cap J = IJ$. We will be using the structure theorem proved by Brumatti and Simis in \cite{bs}; 
see \ref{brusim}. 
\medskip

Let $k$ be a field and $I$ a monomial ideal in the polynomial ring $R=k[x_{1},x_{2},\ldots,x_{n}]$. 
Let us denote its unique minimal generating set by $G(I)$ and $\mbox{Mon}(R)$ denote the set of all monomials 
in $R$. Let us recall some definitions and basic facts from \cite{ssttrans}.

\begin{definition}
Let $\emptyset \neq T \subset$ Mon$(R)$. We define
\begin{center}
supp$(T)=\lbrace i\mid x_{i}$ divides $m$ for some $m\in T\rbrace$.
\end{center}
If $T = \lbrace m\rbrace$, we simply write $\mbox{supp}(m)$ instead of $\mbox{supp}(\lbrace m\rbrace)$. 
Given a nonzero polynomial $f$ in $R$, the set $\mbox{Mon}(f)$ denotes the set of monomials appearing in 
$f$ with nonzero coefficients and $\mbox{supp}(f)$ is defined similarly. If $S$ and $T$
are two nonempty subsets of Mon$(R)$, then, $\mbox{supp}(S)\cap \mbox{supp}(T) = \emptyset$ if and
only if $\mbox{supp}(f) \cap \mbox{supp}(g) = \emptyset$, for every $f\in S$ and $g\in T$.
\end{definition}
\medskip

\begin{theorem}\label{disjoint}
Let $I$ and $J$ be two monomial ideals of $R$. Then, $I \cap J = IJ$
if and only if supp$(G(I))\cap$ supp$(G(J)) = \emptyset$.
\end{theorem}

\proof See theorem $2.2$ in \cite{ssttrans}.\qed
\medskip

\begin{prop}\label{colon}
Let $I$ be a monomial ideal and $G(I)=\lbrace M_{1},\ldots,M_{k}\rbrace$ be the minimal generating set. 
Suppose $m\in R$ be a monomial and $m_{i}=\dfrac{M_{i}}{\gcd(M_{i},m)}$. Then 
$(I:m)=\langle\lbrace m_{1},\ldots,m_{k}\rbrace\rangle$.
\end{prop}

\proof Note that $\gcd(M_{i},m)\mid m$ and therefore $\gcd(M_{i},m)m_{i}\mid mm_{i}$ for all 
$i\in \lbrace 1,\ldots,k\rbrace$. 
Therefore $M_{i}\mid mm_{i}$, i.e., $m_{i}\in (I:m)$ and hence 
$\langle\lbrace m_{1},\ldots,m_{k}\rbrace\rangle\subset (I:m)$. We pick an element $N\in (I:m)$; 
without loss of generality we may assume that $N$ is a monomial. Since $mN\in I$, there exists 
$M_{i}\in G(I)$ such that $M_{i}\mid mN$. Which implies that 
$\dfrac{M_{i}}{\gcd(M_{i},m)} \mid \dfrac{mN}{\gcd(M_{i},m)}$, i.e., $m_{i}\mid \dfrac{mN}{\gcd(M_{i},m)}$. 
Since $\gcd(m_{i},m)=1$, we have $m_{i}\mid N$, therefore $N\in \lbrace m_{1},\ldots,m_{k}\rbrace$.\qed
\medskip

\begin{theorem}[Brumatti, Simis; \cite{bs}]\label{brusim} Let $I\subset R=k[x_{1},x_{2},\ldots,x_{n}]$ be an ideal generated by monomials whose exponents are prime to char $k$. 
Then
$${\rm Der}_{k}(R/I,R/I)=\displaystyle\bigoplus_{i=1}^{n}((I:(I:x_{i}))/I)\dfrac{\partial}{\partial x_{i}} \subset \displaystyle\bigoplus_{i=1}^{n}(R/I)\dfrac{\partial}{\partial x_{i}}.$$
\end{theorem}

\proof See Theorem $2.2.1$ in \cite{bs}.\qed
\medskip

\begin{lemma}\label{der}
Let $I,J\subset R$ be monomial ideals in $R$, with $I\cap J=IJ$. Let $G(I)=\lbrace M_{1},\ldots,M_{k}\rbrace$ 
and $G(J)=\lbrace N_{1},\ldots,N_{s}\rbrace$ be the minimal generating sets of $I$ and $J$ respectively. 
Then
\begin{eqnarray*}
((I+J):((I+J):x_{i})) & =  & 
\begin{cases}
(I:(I:x_{i}))+J \quad \mathrm{if} \quad i\in \mathrm{supp}(G(I))\\
I+(J:(J:x_{i})) \quad \mathrm{if} \quad i\in \mathrm{supp}(G(J))\\
R \quad \mathrm{if} \quad i\notin \mathrm{supp}(G(I))\cup\mathrm{supp}(G(J)).
\end{cases}
\end{eqnarray*}
\end{lemma}

\proof We have $I\cap J=IJ$ if and only if supp$(G(I))\cap$ supp$(G(J))=\emptyset$, 
by theorem \ref{disjoint}. Suppose that $i\in$ $\mbox{supp}(G(I))$, we have 
$(I:x_{i})=\langle m_{i1},\ldots,m_{ik}\rangle$, where $m_{ij}=\dfrac{M_{j}}{\gcd(M_{j},x_{i})}$ by proposition \ref{colon}. 
\medskip

\noindent\textit{Claim 1.} $((I+J):x_{i})=\langle m_{i1},\ldots,m_{ik},N_{1},\ldots,N_{s}\rangle$.
\medskip

\noindent\textit{Proof of Claim 1.} \, It is obvious that $\langle m_{i1},\ldots,m_{ik},N_{1},\ldots,N_{s}\rangle\subset (I+J:x_{i})$. Conversely, 
let $N$ be a monomial such that $x_{i}N\in I+J$. There exists $T\in G(I)\cup G(J)$ such that $T\mid x_{i}N$. 
If $T=M_{j}$ for some $J\in \lbrace 1,\ldots,k\rbrace$, then $M_{j}\mid x_{i}N$ implies that 
$\dfrac{M_{j}}{\gcd(M_{j},x_{i})} \mid \dfrac{x_{i}N}{\gcd(M_{j},x_{i})}$, i.e., $m_{ij}\mid N$. 
If $T=N_{j}$ for some $j\in \lbrace 1,\ldots, s\rbrace$, then $N_{j}\mid x_{i}N$. Therefore $N_{j}\mid N$, 
since $x_{i}\nmid N_{j}$. \qed
\medskip

Now $(I+J:\langle m_{i1},\ldots,m_{ik},N_{1},\ldots,N_{s}\rangle)= (I+J:m_{i1})\cap \cdots\cap (I+J:m_{ik})\cap (I+J:N_{1})\cap \cdots\cap (I+J:N_{s})$. 
\medskip

\noindent\textit{Claim 2.} $\displaystyle \bigcap_{j=1}^{k}((I+J):m_{ij})=(I:(I:x_{i}))+J$.
\medskip

\noindent\textit{Proof of Claim 2.} \, We use Proposition \ref{colon} and get
\begin{align*}
\bigcap_{j=1}^{k}((I+J):m_{ij})&=\bigcap_{j=1}^{k}\left\langle M_{1},\ldots,M_{k},N_{1},\ldots,N_{s}:m_{ij}\right\rangle\\
&=\bigcap_{j=1}^{k}\left\langle\dfrac{M_{1}}{\gcd(M_{1},m_{ij})} ,\ldots,\dfrac{M_{k}}{\gcd(M_{k},m_{ij})},N_{1},\ldots,N_{s}\right\rangle\\
&=\bigcap_{j=1}^{k}\left\langle\dfrac{M_{1}}{\gcd(M_{1},m_{ij})} ,\ldots,\dfrac{M_{k}}{\gcd(M_{k},m_{ij})}\right\rangle +J 
\\
&=(I:(I:x_{i}))+J\qed
\end{align*}

\noindent Therefore, if $ i\in \mathrm{supp}(G(I))$, we have 
$$((I+J):((I+J):x_{i}))=(I:(I:x_{i}))+J.$$ 
Similarly, 
if $i\in \mbox{supp}(G(J))$, we have 
$$((I+J):((I+J):x_{i}))=I+(J:(J:x_{i})).$$ 
If $i\notin\mbox{supp}(G(I))\cup \mbox{supp}(G(J))$, then, by Proposition \ref{colon} 
we get 
$$(I+J:x_{i})=I+J.$$ 
Therefore, $((I+J):((I+J):x_{i}))=R$; this completes the proof of the Lemma. \qed
\medskip

\begin{corollary}\label{directsum}
Let $I,J\subset R=k[x_{1},x_{2},\ldots,x_{n}]$ be ideals generated by monomials whose exponents are prime 
to $\mbox{char}\,k$. Then
$${\rm Der}_{k}(R/I+J,R/I+J)= {\rm Der}_{k}(R/I,R/I)\oplus {\rm Der}_{k}(R/J,R/J).$$
\end{corollary}

\proof Follows from Theorem \ref{disjoint}, Theorem \ref{brusim} and Lemma \ref{der}.\qed
\bigskip

\subsection*{Join of Simplical Complexes}

A \textit{simplicial complex} on $\{1,\ldots,n\}$ is a collection $\Delta$ of subsets of 
$\{1,\ldots,n\}$, such that if $F\in \Delta$ and $F^{'}\subset F$, then $F^{'}\in\Delta$. 
The set $\{1,\ldots,n\}$ is usually called the \textit{vertex set} of $\Delta$, denoted 
by $V(\Delta)$. An element of $\Delta$ is called a face of $\Delta$. A \textit{facet} is a maximal 
face of $\Delta$, with respect to inclusion. Let $\mathcal{F}(\Delta)$ denote the set 
of all facets of $\Delta$. A \textit{nonface} of $\Delta$ is a subset $F$ of $\{1,\ldots,n\}$, 
such that $F\notin \Delta$. Let $\mathcal{N}(\Delta)$ denote the set of minimal nonfaces of 
$\Delta$. For each subset $F\subset \{1,\ldots,n\}$, we set $x_{F}=\displaystyle \prod_{i\in F}x_{i}$. 
The \textit{Stanley-Reisner ideal} of $\Delta$ is the ideal of $R$ defined as 
$I_{\Delta}=\langle \{x_{F}\mid F \in\mathcal{N}(\Delta)\} \rangle$ and the quotient 
ring $k[\Delta]:= R/I_{\Delta}$ is called the \textit{Stanley-Reisner ring} of $\Delta$. 
\medskip

\begin{definition}
Let $\Gamma$ and $\Delta$ be simplicial complexes on disjoint vertex sets $V$ and $W$, 
respectively. The join $\Gamma * \Delta$ is the simpilicial complex on the vertex set $V\cup W$ with 
$\mathcal{F}(\Gamma * \Delta)= \{F\cup G \mid F\in \Gamma, G\in \Delta\} $.
\end{definition}
\medskip

\begin{corollary}\label{dersimplicial}
Let $k$ be a field of characteristic zero. Then 
$${\rm Der}_{k}(k[\Gamma * \Delta],k[\Gamma * \Delta])={\rm Der}_{k}(k[\Delta],k[\Delta])\oplus {\rm Der}_{k}(k[\Gamma],k[\Gamma]).$$
\end{corollary}

\proof Here we note that $I_{\Gamma *\Delta}=I_{\Gamma}+I_{\Delta}$ and $I_{\Gamma}\cap I_{\Delta}=I_{\Gamma}I_{\Delta}$, 
therefore result follows from Corollary \ref{directsum}.\qed
\medskip

\begin{example}
Let $\Gamma$ be a simplicial complex with $V(\Gamma)=\{1,2,3,4,5\}$ and 
$\mathcal{F}(\Gamma)=\{\{1,2,4\},\{1,2,5\},\{2,3\},\{3,4\}\}$. 
Let $\Delta$ denote the simplicial complex with $V(\Delta)=\{6,7,8,9,10\}$ 
and $\mathcal{F}(\Delta)=\{ \{6,7,8\},\{7,9\},\{8,9\},\{10\}\}$. 
One can show 
after some easy computation that 
$$D(k[\Gamma * \Delta],k[\Gamma * \Delta])=D(k[\Delta],k[\Delta])\oplus D(k[\Gamma],k[\Gamma]).$$
\end{example}
\bigskip

\section{Derivation modules of glued semigroup rings}

A \textit{numerical semigroup} $\Gamma$ is a subset of the set of nonnegative integers 
$\mathbb{N}$, closed under addition, contains zero and generates $\mathbb{Z}$ as 
a group. It follows that (see \cite{rosp}) the set $\mathbb{N}\setminus \Gamma$ is 
finite and that the semigroup $\Gamma$ has a unique minimal system of generators 
$n_{0} < n_{1} < \cdots < n_{p}$. The greatest integer not belonging to $\Gamma$ 
is the \textit{Frobenius number}, $n_{0}$ is the \textit{multiplicity} 
and $p + 1$ is the \textit{embedding dimension} of the numerical semigroup $\Gamma$. 
The \textit{Ap\'{e}ry set} of $\Gamma$ with respect to a non-zero $a\in \Gamma$ is 
the set $\rm{Ap}(\Gamma,a)=\{s\in \Gamma\mid s-a\notin \Gamma\}$. 
\medskip

Let $p\geq 1$ and $n_{0}, \ldots, n_{p}$ be positive integers with 
$\gcd (n_{0},\ldots,\, n_{p})=1$. Let us assume that no $n_{i}$ can 
be written in terms of linear combination of other $n_{j}$'s over 
$\mathbb{Z}_{\geq 0}$. The set 
$\Gamma(n_{0},\ldots, n_{p}):= \{\sum_{i=0}^{p}c_{i}n_{i} \mid c_{i}\in\mathbb{Z}_{\geq 0}\}$ 
is a numerical semigroup and it is minimally generated by $n_{0}, \ldots, n_{p}$. 
Let $k$ denote a field and $R$ denote the polynomial ring $k[x_0,\,\ldots,\, x_p]$. 
One can define a $k$-algebra homomorphism $\eta:k[x_{0},\ldots, x_{p}]\rightarrow k[t]$ as 
$\eta(x_i)=t^{n_i},\,0\leq i\leq p$. Let $\frak{p}(n_0,\ldots, n_p) = 
\ker (\eta)$. The map $\eta$ is a parametrization of a curve, 
known as a \textit{monomial curve} in the affine space $\mathbb{A}^{p+1}$. 
The ideal $\frak{p}(n_0,\ldots, n_p)$ is called the defining ideal of the 
curve $\mathcal{C}$. The affine $k$-algebra 
$R_{\Gamma} = k[x_0,\,\ldots,\, x_p]/\frak{p}(n_0,\ldots, n_p)=k[t^{n_{0}}, \ldots ,t^{n_{p}}]$ 
is called the \textit{coordinate ring} of the curve $\mathcal{C}$. 
We fix some notations first: 
\medskip

\begin{itemize}
\item $\Gamma(n_0,\ldots, n_p)_{+} := \Gamma(n_0,\ldots, n_p) \setminus \{0\};$
\item $\Delta_{\Gamma} := \{ \alpha \in \mathbb{Z}^{+} \mid \alpha + \Gamma(n_0,\ldots, n_p)_{+} \subseteq    
\Gamma(n_0,\ldots, n_p)\}; $
\item $\Delta_{\Gamma}^{'} := \Delta_{\Gamma}\setminus \Gamma(n_0,\ldots, n_p).$
\end{itemize}
\medskip

\begin{lemma}\label{delta}
Let $\Gamma$ be a numerical semigroup with $a\in \Delta_{\Gamma}^{'}$ and $b\in \Gamma$, then, 
$\Delta_{\Gamma}^{'}\subset\mathrm{Ap}(\Gamma,b)-b$.
\end{lemma}

\proof There exists $w\in \mathrm{Ap}(\Gamma,b)$, such that $a\equiv w(\mathrm{mod}\, b)$. If 
$a\geq w$, then $a-w=kb$ for some $k\geq 0$. This implies that $a=w+kb\in \Gamma$, which is a 
contradiction. Therefore, $a<w$. Let $w-a=kb$, where $k>0$ (since $a\equiv w(\mathrm{mod}\, b))$. 
If $k>1$, then $w-b=(k-1)b+a$. Since $a\in \Delta_{\Gamma}^{'}$ and $(k-1)>0$, therefore 
$a+(k-1)b\in \Gamma$ implies that $w-b\in \Gamma$, which gives a contradiction to the fact that 
$w\in \mathrm{Ap}(\Gamma,b)$. Hence $k=1$, i.e., $a=w-b$. \qed
\medskip

\begin{definition}\label{gluedef} Let $\Gamma_{1}$ and $\Gamma_{2}$ be two numerical 
semigroups minimally generated by $\{n_{1},\ldots,n_{r}\}$ and $\{n_{r+1},\ldots,n_{e}\}$, 
respectively. Let $\lambda\in \Gamma_{1}\setminus\{n_{1},\ldots,n_{r}\}$ and 
$\mu\in \Gamma_{2}\setminus\{n_{r+1},\ldots,n_{e}\}$
be such that $\gcd(\lambda,\mu)=1$. We say that
$\Gamma := \langle \{\mu n_{1},\ldots,\mu n_{r},\lambda n_{r+1},\ldots,\lambda n_{e}\}\rangle$
is a \textit{gluing} of the numerical semigroups $\Gamma_{1}$ and $\Gamma_{2}$ and it is 
denoted by $\mu \Gamma_{1}\# \lambda \Gamma_{2}$.
\end{definition}
\medskip

\begin{lemma}\label{glueapery}
Let $\Gamma_{1}$ and $\Gamma_{2}$ be two numerical 
semigroups minimally generated by $\{n_{1},\ldots,n_{r}\}$ and $\{n_{r+1},\ldots,n_{e}\}$, 
respectively. Let $\lambda\in \Gamma_{1}\setminus\{n_{1},\ldots,n_{r}\}$ and 
$\mu\in \Gamma_{2}\setminus\{n_{r+1},\ldots,n_{e}\}$
be such that $\gcd(\lambda,\mu)=1$.
\begin{enumerate}
\item The numerical semigroup $\Gamma = \mu \Gamma_{1}\# \lambda \Gamma_{2}$ is minimally generated by set 
$\{\mu n_{1},\ldots,\mu n_{r},\lambda n_{r+1},\ldots,\lambda n_{e}\}$.

\item $\mathrm{Ap}(\Gamma ,\lambda\mu)=\{\mu w_{1}+\lambda w_{2}\mid w_{1}\in 
\mathrm{Ap}(\Gamma_{1},\lambda),\, w_{2}\in \mathrm{Ap}(\Gamma_{2},\mu)\}$ 
\end{enumerate}
\end{lemma}

\proof See Lemma 9.8, Proposition 9.11 and Theorem 9.2 in \cite{rosp}.  \qed
\medskip 

Let ${\rm Der}_{k}(R_{\Gamma})$ denote the set of all $k$-derivations of $R_{\Gamma}$, 
called the \textit{Derivation module} of $R_{\Gamma}$. It is a finitely generated 
module over $R_{\Gamma}$. The following theorem by J. Kraft in \cite{kraft}, 
page 875, gives an explicit minimal generating set of the $R_{\Gamma}$-module 
${\rm Der}_{k}(R_{\Gamma})$. 
\medskip 

\begin{theorem}[Kraft; \cite{kraft}] \label{kraft} The set $\left\{ t^{\alpha + 1} \frac{d}{dt} \mid \alpha \in
\Delta'_{\Gamma} \cup \{0\}\right\}$ is a minimal set of generators 
for the $R_{\Gamma}$-module ${\rm Der}_{k}(R_{\Gamma})$, where $\Gamma(n_0,\ldots, n_p)$ is 
the numerical semigroup minimally generated by the sequence $n_0,\ldots, n_p$ 
of positive integers. In particular, $\mu ({\rm Der}_{k}(R_{\Gamma})) = {\rm card}\left(\Delta_{\Gamma}^{'}\right) +1$. 
\end{theorem}

\proof See \cite{kraft} page no 875.\qed
\medskip

We now prove the following theorem, which gives the description of $\Delta_{\Gamma}^{'}$, 
where $\Gamma$ is the numerical semigroup obtained by gluing the numerical semigroups 
$\Gamma_{1}$ and $\Gamma_{2}$.
\medskip

\begin{theorem}\label{deltaglue}
Let $\Gamma_{1}$ and $\Gamma_{2}$ be two numerical 
semigroups minimally generated by $\{n_{1},\ldots,n_{r}\}$ and $\{n_{r+1},\ldots,n_{e}\}$, 
respectively. Let $\lambda\in \Gamma_{1}\setminus\{n_{1},\ldots,n_{r}\}$ and 
$\mu\in \Gamma_{2}\setminus\{n_{r+1},\ldots,n_{e}\}$
be such that $\gcd(\lambda,\mu)=1$. Let $\Gamma$ be their glued semigroup, i.e., 
$\Gamma = \mu \Gamma_{1}\# \lambda \Gamma_{2}$. If 
$\Delta_{\Gamma_{1}}^{'}=\lbrace a_{1},a_{2},\ldots, a_{n}\rbrace$ and 
$\Delta_{\Gamma_{2}}^{'}=\lbrace b_{1},b_{2},\ldots, b_{m}\rbrace$, then,  
$\Delta_{\Gamma}^{'}=\lbrace a_{i}\mu+b_{j}\lambda+\lambda\mu \mid 1\leq i\leq n, 1\leq j\leq m\rbrace$. 
Hence 
${\rm card}\left(\Delta_{\Gamma}^{'}\right)={\rm card}\left(\Delta_{\Gamma_{1}}^{'}\right)\cdot{\rm card}\left(\Delta_{\Gamma_{2}}^{'}\right)$.
\end{theorem}

\proof It is enough to prove the statements (i), (ii), (iii) and (iv) given below:

\begin{itemize}
\item[(i)] For every $1\leq i\leq n$, $1\leq j\leq m$, we have 
$\mu a_{i}+\lambda b_{j}+\lambda\mu \notin \Gamma$.
\medskip

\item[(ii)] For $t\in \Gamma\setminus\lbrace0\rbrace$, $\mu a_{i}+\lambda b_{j}+\lambda\mu +t\in \Gamma$.
\medskip

\item[(iii)] Let $a\notin \Gamma$, $a+\mu n_{i}\in \Gamma$ and $a+\lambda n_{j}\in \Gamma$, where 
$1\leq i\leq r$, $r+1\leq j\leq e$. Then, \, $a$ is of the form $\mu a_{i}+\lambda b_{j}+\lambda\mu$;  
$1\leq i\leq n, 1\leq j\leq m$.
\medskip

\item[(iv)] Given $i,i^{\prime}\in \lbrace 1,\ldots, n\rbrace$, \,\, $j,j^{\prime}\in \lbrace 1,\ldots, m\rbrace$ 
and $(i,j)\neq (i^{'},j^{'})$, we have $\mu a_{i}+\lambda b_{j}+\lambda\mu \neq \mu a_{i^{\prime}}+\lambda b_{j^{\prime}}+\lambda\mu$. 
\end{itemize}
\medskip

\noindent\textbf{Proof of (i).} Suppose we consider $\mu a_{1}+\lambda b_{2}+\lambda\mu\in \Gamma$, then, we have 
\begin{align*}
\mu a_{1}+\lambda b_{2}+\lambda\mu &= \displaystyle\sum_{i=1}^{r} d_{i}\mu n_{i}+\displaystyle\sum_{i=r+1}^{e} d_{i}\lambda n_{i}.
\end{align*}
Therefore we get,
\begin{align*}
\mu (a_{1}-\displaystyle\sum_{i=1}^{r} d_{i}n_{i})&= \lambda(\displaystyle\sum_{i=r+1}^{e} d_{i}n_{i}-b_{2}-\mu).
\end{align*}
If $a_{1}>\displaystyle\sum_{i=1}^{r} d_{i} n_{i}$, then $\lambda\vert(a_{1}-\displaystyle\sum_{i=1}^{r} d_{i}n_{i})$
and $a_{1}-\displaystyle\sum_{i=1}^{r} d_{i}n_{i}=\lambda k$, $k>0$. Therefore,
$ a_{1}=\displaystyle\sum_{i=1}^{r} d_{i}n_{i}+\lambda k\in \Gamma_{1}$, since 
$\displaystyle\sum_{i=1}^{r} d_{i}n_{i}\in \Gamma_{1}$ and $\lambda k\in \Gamma_{1}$; which is a contradiction.
\medskip

Next we consider the case \, $a_{1}<\displaystyle\sum_{i=1}^{r} d_{i} n_{i}$.
Then we have, $$\lambda(b_{2}+\mu-\displaystyle\sum_{i=r+1}^{e} d_{i} n_{i})=\mu (\displaystyle\sum_{i=1}^{r} d_{i} n_{i}-a_{1}).$$
Therefore $\mu\vert(b_{2}+\mu-\displaystyle\sum_{i=r+1}^{e} d_{i} n_{i})$, which implies that 
$b_{2}+\mu-\displaystyle\sum_{i=r+1}^{e} d_{i} n_{i}=k^{\prime}\mu$, for some  $k^{\prime}>0$. 
Hence $b_{2}=(k^{\prime}-1)\mu+\displaystyle\sum_{i=r+1}^{e} d_{i} n_{i}$, $k^{\prime}-1\geq 0$, 
which implies that $b_{2}\in S_{2}$, since $(k^{\prime}-1)\mu\in \Gamma_{2}$ and 
$\displaystyle\sum_{i=r+1}^{e} d_{i} n_{i}\in \Gamma_{2}$. This is again a contradiction.
\medskip

\noindent\textbf{Proof of (ii).} $t\in \Gamma\setminus\lbrace 0\rbrace$ implies that 
$t=\displaystyle\sum_{i=1}^{r} d_{i}\mu n_{i}+\displaystyle\sum_{i=r+1}^{e} d_{i}\lambda n_{i}$.
Now we have $\mu a_{i}+\lambda b_{j}+\lambda\mu+t = \mu a_{i}+\lambda b_{j}+\lambda\mu+\displaystyle\sum_{i=1}^{r} d_{i}\mu n_{i}+\displaystyle\sum_{i=r+1}^{e} d_{i}\lambda n_{i}$. Therefore
$\mu a_{i}+\lambda b_{j}+\lambda\mu+t = \mu(a_{i}+\displaystyle\sum_{i=1}^{r} d_{i}n_{i})+\lambda(b_{j}+\displaystyle\sum_{i=r+1}^{e} d_{i}n_{i})+\lambda\mu \in \Gamma$.
\medskip

\noindent\textbf{Proof of (iii).}  Let $a\in \Delta_{\Gamma}^{'}$. By lemma \ref{delta} 
there exists $w\in \mathrm{Ap}(\Gamma,\lambda\mu)$ such that $w=a+\lambda\mu$. We know that  
$w\in \mathrm{Ap}(\Gamma,\lambda\mu)$, therefore, by lemma \ref{glueapery} we have 
$w=\mu w_{1}+\lambda w_{2}$, where $w_{1}\in \mathrm{Ap}(\Gamma_{1},\lambda)$ and 
$w_{2}\in \mathrm{Ap}(\Gamma_{2},\mu)$. Therefore 
$a=\mu w_{1}-\lambda\mu+\lambda w_{2}=\mu(w_{1}-\lambda)+\lambda w_{2}$. Note that 
$w_{1}\in \mathrm{Ap}(\Gamma_{1},\lambda)$ implies $w_{1}-\lambda\notin \Gamma_{1}$. 
Now we have,
\begin{align*}
a+\mu n_{i}=\mu(w_{1}-\lambda+n_{i})+\lambda w_{2},\, \mathrm{for\, all}\, 1\leq i\leq r .
\end{align*}
We claim that, $w_{1}-\lambda+n_{i}\in \Gamma_{1}$, for all $1\leq i\leq r$. We know that 
$w_{1}+n_{i}\in \Gamma_{1}$ for all $1\leq i\leq r$, therefore, it is enough to show that 
$w_{1}+n_{i}\notin \mathrm{Ap}(\Gamma_{1},\lambda)$, for all $1\leq i\leq r$. 
\medskip

We have $a+\mu n_{i}= \mu (w_{1}-\lambda+n_{i})+\lambda w_{2}\in \Gamma$, for all $1\leq i\leq r$. 
Suppose that $w_{1}+n_{i}\in \mathrm{Ap}(\Gamma_{1},\lambda)$. Then 
$\mu(w_{1}+n_{i})+\lambda w_{2}\in \mathrm{Ap}(\Gamma,\lambda\mu)$ by lemma\ref{glueapery}, 
which implies that $\mu(w_{1}+n_{i})+\lambda w_{2}-\lambda\mu\notin S$, i.e., 
$\mu(w_{1}-\lambda+n_{i})+\lambda w_{2}\notin \Gamma$, which is a contradiction. Now we have 
$w_{1}-\lambda\notin \Gamma_{1}$ and $w_{1}-\lambda+n_{i}\in \Gamma_{1}$ for all $1\leq\ i\leq\ r$, 
so that $w_{1}-\lambda=a_{i}$ for some $i\in \lbrace 1,\ldots,n\rbrace$. 
Therefore, $a=\mu a_{i}+\lambda w_{2}= \mu a_{i}+\lambda(w_{2}-\mu)+\lambda\mu$. Similarly we get $w_{2}-\mu =b_{j}$ 
for some $j\in \lbrace 1,\ldots,m\rbrace$. Hence $a=\mu a_{i}+\lambda b_{j}+\lambda\mu$ for some $i\in \lbrace 1,\ldots,n\rbrace$ 
and $j\in \lbrace 1,\ldots,m\rbrace$.
\medskip

\noindent\textbf{Proof of (iv).} Without loss of generality we may assume that 
$a_{i}>a_{i^{\prime}}$. Suppose that \, $\mu a_{i}+\lambda b_{j}+\lambda\mu = \mu a_{i^{\prime}}+\lambda b_{j^{\prime}}+\lambda\mu$, 
then we have $ \mu(a_{i}-a_{i^{\prime}})=\lambda(b_{j^{\prime}}-b_{j})$. Since$(\mu,\lambda)=1$, 
we have $\lambda\vert(a_{i}-a_{i^{\prime}})$ and $\mu\vert(b_{j^{\prime}}-b_{j})$. Let $a_{i}-a_{i^{\prime}}= k\lambda$, 
for some $k\geq 0$. If $k=0$, then $a_{i}=a_{i^{\prime}}$, which is a contradiction. If $k>0$, then 
$a_{i}=a_{i^{\prime}}+k\lambda\in \Gamma_{1}$. Now $k\lambda\in \Gamma_{1}$ and 
$a_{i^{\prime}}\in \Delta_{\Gamma_{1}}^{\prime}$ imply that $a_{i}\in \Gamma_{1}$, which is absurd 
since $a_{i}\in\Delta_{\Gamma_{1}}^{\prime}$. \qed
\medskip

\begin{example}
Let $\Gamma_{1}=\langle 4,7,9\rangle$ and $\Gamma_{2}=\langle 3,8\rangle$ be two numerical semigroups. 
Take $\lambda =8\in S_{1}\setminus\lbrace 4,7,9\rbrace$ and $\mu =11\in S_{2}\setminus\lbrace 3,8\rbrace$ 
so that $\gcd(\lambda,\mu)=1$. The glued semigroup 
$\Gamma=\mu \Gamma_{1}\#\lambda \Gamma_{2}=\langle 44,77,99,24,64\rangle$. It not difficult to show that 
$\Delta_{\Gamma_{1}}^{'}=\lbrace 5,10\rbrace, \Delta_{\Gamma_{2}}^{'}=\lbrace 13\rbrace$ and 
$\Delta_{\Gamma}^{'}=\lbrace 247,302\rbrace$. We observe that $247=5\times 11+13\times 8+8\times 11$ and $302=10\times 11+13\times 8+8\times 11$.
\end{example}
\medskip

\begin{corollary}
Let $\Gamma_{1}$ and $\Gamma_{2}$ be two numerical semigroups minimally generated by 
$\{n_{1},\ldots,n_{r}\}$
and $\{n_{r+1},\ldots,n_{e}\}$, respectively. Let $\lambda\in \Gamma_{1}\setminus\{n_{1},\ldots,n_{r}\}$ 
and $\mu\in \Gamma_{2}\setminus\{n_{r+1},\ldots,n_{e}\}$
be such that $\gcd(\lambda,\mu)=1$. Suppose $\Gamma$ be their glued semigroup, i.e., 
$\Gamma=\mu \Gamma_{1}\# \lambda \Gamma_{2}$. Then the set 
$\left\{ t^{a\mu+b\lambda+\lambda\mu+ 1} \frac{d}{dt} \mid a \in
\Delta'_{\Gamma_{1}}, b \in
\Delta'_{\Gamma_{2}} \right\}\cup \{t\frac{d}{dt}\}$ is a minimal set of generators 
for the $R_{\Gamma}$-module ${\rm Der}_{k}(R_{\Gamma})$ and 
$$\mu({\rm Der}_{k}(R_{\Gamma})) = {\rm card}\left(\Delta_{\Gamma}^{'}\right) +1 
= {\rm card}\left(\Delta_{\Gamma_{1}}^{'}\right)\cdot{\rm card}\left(\Delta_{\Gamma_{2}}^{'}\right) + 1.$$
\end{corollary}

\proof Follows from theorems \ref{deltaglue} and \ref{kraft}.\qed
\medskip

\begin{theorem}\label{ghglue} 
Suppose that a numerical semigroup $\langle \Gamma \rangle$ is obtained by gluing the numerical semigroups 
$\Gamma_{1}$ and $\Gamma_{2}$. 
Let $\beta(R_{\Gamma})$, $\beta(R_{\Gamma_{1}})$ and $\beta(R_{\Gamma_{2}})$ denote the last Betti numbers of the 
numerical semigroup rings $R_{\Gamma}$, $R_{\Gamma_{1}}$ and $R_{\Gamma_{2}}$ respectively. Then 
$$\beta(\Gamma)={\rm card}\left(\Delta_{\Gamma_{1}}^{'}\right)\cdot{\rm card}\left(\Delta_{\Gamma_{2}}^{'}\right) = \beta(\Gamma_{1}).\beta(\Gamma_{2}).$$
\end{theorem}

\proof For every numerical semigroup $\Gamma$, it is true that 
$\beta(R_{\Gamma})=t_{R_{\Gamma}}=\mu({\rm Der}_{k}(R_{\Gamma}))-1$; see corollary 6.2 in \cite{ps} for a proof. Therefore, 
\begin{eqnarray*}
\beta(R_{S}) & = & \mu({\rm Der}_{k}(R_{\Gamma}))-1\\
{} & = & {\rm card}\left(\Delta_{\Gamma_{1}}^{'}\right)\cdot{\rm card}\left(\Delta_{\Gamma_{2}}^{'}\right)\\
{} & = & (\mu({\rm Der}_{k}(R_{\Gamma_{1}}))-1)\cdot (\mu({\rm Der}_{k}(R_{\Gamma_{2}}))-1)\\
{} & = & \beta(R_{\Gamma_{1}}).\beta(R_{\Gamma_{2}}). 
\end{eqnarray*}
This result also follows from Theorem 3.1 in \cite{gs} proved by Gimenez and Srinivasan. However, an explicit 
computation of $\beta(R_{\Gamma})$ is possible from their result if $\beta(R_{\Gamma_{1}})$ and 
$\beta(R_{\Gamma_{2}})$ 
are known and that is not easy to compute through homological tools. On the other hand, we only need to 
calculate ${\rm card}\left(\Delta_{\Gamma_{1}}^{'}\right)$ and ${\rm card}\left(\Delta_{\Gamma_{2}}^{'}\right)$ 
in order to calculate $\beta(R_{\Gamma})$. \qed
\bigskip

\noindent\textbf{The GAP program\footnote{The authors thank Amogh Parab 
for writing this program, which has helped immensely in calculating $\Delta_{S}^{'}$ for various test 
cases at the formative stage of this work.}
for calculating $\Delta_{S}^{'}$}
\medskip

\noindent $\mathrm{gap>\,p:=[n_{1},\ldots,n_{e}];\quad\quad // \lbrace n_{1},\ldots,n_{e}\rbrace \,are\, minimal\, generators\, of\, S}$\\
$\mathrm{gap> m:=n_{1}};$\\
$\mathrm{gap>S:=NumericalSemigroup(p)};$\\
$\mathrm{gap>AP:=AperyListOfNumericalSemigroupWRTElement(S,m);}$\\
$\mathrm{gap>AP[1]:=m;}$\\
$\mathrm{gap>F:=FrobeniusNumberOfNumericalSemigroup(S);}$\\
$\mathrm{gap>L:=[\,];}$\\
$\mathrm{gap>a:=0;}$\\
$\mathrm{gap>while\, a<F+1\, do}$\\
$\mathrm{>if\, a\, in\, S\, then}$\\
$\mathrm{>d:=1;}$\\
$\mathrm{>else}$\\
$\mathrm{>c:=0;}$\\
$\mathrm{>for\, i\, in\, [1..m]\, do}$\\
$\mathrm{>if\, a+AP[i]\, in\, S\, then}$\\
$\mathrm{>c:=c+1;}$\\
$\mathrm{>fi;}$\\
$\mathrm{>od;}$\\
$\mathrm{>if\, c=m\, then}$\\
$\mathrm{>Append(L, [a]);}$\\
$\mathrm{>fi;}$\\
$\mathrm{>fi;}$\\
$\mathrm{>a:=a+1;}$\\
$\mathrm{>od;}$\\
$\mathrm{gap>L; \quad\quad // L\, will\, give\, \Delta_{S}^{'}}$

\bibliographystyle{amsalpha}

\end{document}